\documentclass[12pt,a4wide]{article}
\usepackage[leqno]{amsmath}
\usepackage[all]{xy}
\usepackage[mathscr]{euscript}
\usepackage{amssymb}
\usepackage{amscd}

\newtheorem{lemma}{Lemma}[section]
\newtheorem{theorem}[lemma]{Theorem}
\newtheorem{prop}[lemma]{Proposition}

\newtheorem{cor}[lemma]{Corollary}

\newcommand{\pf}{\noindent{\em Proof: }}
\newcommand{\epf}{\hfill\hbox{\rule{3pt}{6pt}}\\}

\newcommand{\forme}[1]{}
\begin{document}

\title{An inequality involving the second largest and smallest eigenvalue of a distance-regular graph }

\author{Jack H. Koolen$^{\,\rm 1,2,}$\footnote{
This work was partially supported by the Priority Research Centers Program through the National Research Foundation of Korea (NRF) funded by the Ministry of Education, Science and Technology (Grant number 2009-0094069). JHK was also partially supported by the Basic Science Research Program through the National Research
Foundation of Korea(NRF) funded by the Ministry of Education, Science and Technology (Grant number 2009-0089826).
}
,~  Jongyook Park$^{\,\rm 1}$
and Hyonju Yu$^{\,\rm 1}$\\
{\small {\tt koolen@postech.ac.kr} ~~
{\tt jongyook@postech.ac.kr} ~~ {\tt lojs4110@postech.ac.kr}}\\
{\footnotesize{$^{\rm 1}$Department of Mathematics,  POSTECH, Pohang 790-785, South Korea}}\\
{\footnotesize{$^{\rm 2}$Pohang Mathematics Institute,  POSTECH, Pohang 790-785, South Korea}}}

\date{\today}

\maketitle

\begin{abstract} For a distance-regular graph with second largest eigenvalue (resp. smallest eigenvalue) $\theta_1 \ $ (resp. $\theta_D$) we show that
$ (\theta_1 + 1)(\theta_D +1) \leq -b_1$ holds, where equality only holds when the diameter equals two. Using this inequality we study distance-regular graphs with fixed second largest eigenvalue.

\end{abstract}

\section{Introduction}
In Juri\v{s}i\'{c} et al. \cite{JKT00}, it was shown that for a distance-regular graph with diameter $D$ at least two  and distinct eigenvalues $k = \theta_0 > \theta_1 > \ldots > \theta_D$, one has the following bound:
$$\Biggl( \theta_1 + \frac{k}{a_1 + 1}\Biggr) \Biggl( \theta_D + \frac{k}{a_1 +1}\Biggr) \geq \frac{ka_1b_1}{(a_1 +1)^2}.$$
(For definitions see Section 2.)
\\
\\
In this paper we show another bound involving the second largest and smallest eigenvalue, namely we show:
$ ( \theta_1 + 1)(\theta_D +1) \leq -b_1,$ where equality  holds when and only when the diameter is equal to two (Theorem 3.6).

In order to show this new bound, we give some bounds on the second largest eigenvalue and smallest eigenvalue, respectively, in Section 3.

In Section 4, we determine the distance-regular graphs with diameter at least three and second largest eigenvalue at most two (Theorems 4.1 and 4.3). Also we show  that,  for a fixed integer $m$ at least two, there are only finitely many distance-regular graphs with diameter at least three whose second largest eigenvalue lies in the half-open interval $(1, m]$  (Theorem 4.2).
\section{Definitions and preliminaries}
All the graphs considered in this paper are finite, undirected and
simple (for unexplained terminology and more details, see \cite{bcn}). Suppose
that $\Gamma$ is a connected graph with vertex set $V(\Gamma)$ and edge set $E(\Gamma)$, where $E(\Gamma)$ consists of unordered pairs of two adjacent vertices. The distance $d(x,y)$ between
any two vertices $x,y$ of $\Gamma$
is the length of a shortest path connecting $x$ and $y$ in $\Gamma$.

Let $\Gamma$ be a connected graph. For a vertex $x \in V(\Gamma)$, define $\Gamma_i(x)$ to be the set of
vertices which are at distance precisely $i$ from $x~(0\le i\le
D)$, where $D:=\max\{d(x,y)\mid x,y\in V(\Gamma)\}$ is the diameter
of $\Gamma$. In addition, define $\Gamma_{-1}(x) = \Gamma_{D+1}(x)
= \emptyset$. We write $\Gamma(x)$ instead of $ \Gamma_1(x)$ and we denote $x\sim_{\Gamma} y$ or simply $x\sim y$ if two vertices $x$ and $y$ are adjacent in $\Gamma$. The adjacency matrix $A$ of graph $\Gamma$ is the (0,1)-matrix whose rows and columns are indexed by the vertex set $V(\Gamma)$ and the $(x,y)$-entry is $1$ whenever $x\sim y$ and 0 otherwise.
The eigenvalues (respectively the spectrum) of the graph $\Gamma$ are the eigenvalues (respectively the spectrum) of  $A$.

For a connected graph $\Gamma$, the {\em local graph} $\Delta(x)$ at a vertex $x\in V(\Gamma)$ is the subgraph induced on $\Gamma(x)$ in $\Gamma$.

For a graph $\Gamma$, a partition $\Pi= \{P_1, P_2, \ldots ,
P_{\ell}\}$ of the vertex set $V(\Gamma)$ is called {\em
equitable} if there are constants $\beta_{ij}$ such that each vertex $x \in P_i$ has exactly $\beta_{ij}$ neighbors in $P_j$ ($1\leq i, j \leq \ell$). The {\em quotient
matrix} $Q(\Pi)$ associated with the equitable partition $\Pi$ is
the $\ell \times \ell$ matrix whose $(i,j)$-entry $Q(\Pi)_{(i,j)}$ is $\beta_{ij}$ ($1\leq i,j\leq \ell$). Note that the eigenvalues of the quotient matrix $Q(\Pi)$ are also eigenvalues (of the adjacency matrix $A$) of $\Gamma$ \cite[Theorem 9.3.3]{godsil-combin}.

A connected graph $\Gamma$ with diameter $D$ is called {\em{distance-regular}} if there are integers $b_i, c_{i+1}$ $(0\leq i\leq D-1)$ such that for any two vertices $x, y \in V(\Gamma)$ with $d(x, y)=i$, there are precisely $c_i$ neighbors of $y$ in $\Gamma_{i-1}(x)$ and $b_i$ neighbors of $y$ in $\Gamma_{i+1}(x)$. In particular, any distance-regular graph  is regular with valency $k := b_0$. Note that a {\em non-complete, connected  strongly regular graph} is just a distance-regular graph with diameter two.
We define $a_i := k-b_i-c_i$ for notational convenience.  Note that $a_i=\mid\Gamma(y)\cap\Gamma_i(x)\mid$ holds for any two vertices $x, y$ with $d(x, y)=i$ $(0\leq i\leq D).$ For a distance-regular graph $\Gamma$ and a vertex $x\in V(\Gamma)$, we denote $k_i:=|\Gamma_i(x)|$. It is easy to see that $k_i = \frac{b_0 b_1 \cdots b_{i-1}}{c_1 c_2 \cdots c_i}$ and hence  $k_i$ does not depend on the vertex $x$.  The numbers $a_i$, $b_{i-1}$ and $c_i$ $(1\leq i\leq D)$ are called the {\em{intersection~numbers}} of the distance-regular graph $\Gamma$, and the array $\{b_0,b_1,\ldots,b_{D-1};c_1,c_2,\ldots,c_D\}$ is called the {\em{intersection~array}} of $\Gamma$.

Some easy properties of the intersection numbers are collected in the following lemma.
\begin{lemma}(\cite[Proposition 4.1.6]{bcn})\label{pre}{\ \\}
Let $\Gamma$ be a distance-regular graph with valency $k$ and diameter $D$. Then the following holds:\\
(i) $k=b_0> b_1\geq \cdots \geq b_{D-1}~;$\\
(ii) $1=c_1\leq c_2\leq \cdots \leq c_{D}~;$\\
(iii) $b_i\ge c_j$ \mbox{ if }$i+j\le D~.$
\end{lemma}

Suppose that $\Gamma$ is a distance-regular graph with valency $k\ge 2$ and diameter $D\ge 1$. Then $\Gamma$ has exactly $D+1$ distinct eigenvalues, $k=\theta_0>\theta_1>\cdots>\theta_D$ (\cite[p.128]{bcn}), and the multiplicity of
$\theta_i$ ($0\le i\le D$) is denoted by $m_i$. For an eigenvalue $\theta$ of $\Gamma$, the sequence
$(u_i)_{i=0,1, \ldots, D} = (u_i(\theta))_{i=0,1 , \ldots, D}$ satisfying $u_0=u_0(\theta)=1, u_1=u_1(\theta)=\theta /k, $ and
\begin{equation}\label{rec-d}
c_i u_{i-1} + a_i u_i + b_i u_{i+1} = \theta u_i \ \ (i=2, 3, \ldots, D-1)
\end{equation}
is called the {\em standard sequence} corresponding to the eigenvalue $\theta$ (\cite[p.128]{bcn}).
A {\em sign change} of $(u_i)_{i= 0,1, \ldots, D}$ is a pair $(i,j)$ with $0 \leq i < j \leq D$ such that $u_i u_j < 0$ and
$u_t = 0$ for $i < t <j$.

For a distance-regular graph $\Gamma$ with diameter $D$ and eigenvalues $\theta_0 > \theta_1 > \cdots > \theta_D$, define the tridiagonal $D \times D$ matrix $T= T(\Gamma)$ by
$$T =\left[
\begin{array}{cccccc}
 -c_1 & b_1 & & & & \\
 c_1 & k- b_1 -c_2 & b_2 & & \bf{0} &  \\
 & c_2 & \cdot & \cdot & & \\
 & & \cdot & \cdot & \cdot & \\
 &\bf{0} &  & \cdot & \cdot & b_{D-1} \\
 &&&& c_{D-1} & k-b_{D-1}-c_D
\end{array} \right].$$
It is known that $\theta_1, \theta_2, \cdots, \theta_D$ are exactly the eigenvalues of the matrix $T$ (see for example \cite[p.130]{bcn}).
\\
\\
Recall the following interlacing result.

\begin{theorem}\label{0}\rm\textbf{(cf. Haemers \cite{heam})}   Let $m \geq n$ be two positive integers. Let $A$ be an  $n\times n$ matrix, which is similar to a (real) symmetric matrix, and let $B$ be a principal $m \times m$ submatrix of $A$. Then, for $i=1,\ldots , m$, $$\theta_{n-m+i}(A)\leq \theta_i(B)\leq \theta_i(A)$$
holds, where $A$ has eigenvalues $\theta_1(A) \geq \theta_2(A) \geq \cdots \geq \theta_n(A)$ and $B$ has eigenvalues  $\theta_1(B) \geq \theta_2(B) \geq \cdots \geq \theta_m(B)$.
\end{theorem}

\section{Some inequalities}
In this section we will give some inequalities for the eigenvalues of a distance-regular graph.
For a distance-regular graph $\Gamma$ and a vertex $x$ of $\Gamma$, the subgraph induced on the set $\Gamma(x) \cup \{ x\}$ is denoted by $\widehat{\Delta(x)}$.
\begin{theorem}\label{ineq1}
Let $\Gamma$ be a distance-regular graph
with valency $k$ at least three, diameter $D$ at least three and distinct eigenvalues $k = \theta_0 > \theta_1 > \cdots > \theta_D$.
Then the following hold:
\\
(i) $\theta_D < \frac{a_1 - \sqrt{a_1^2 + 4k}}{2}$;
\\
(ii) $\theta_1 \geq \min\{  \frac{a_1 + \sqrt{a_1^2 + 4k}}{2}, a_3\}$;
\\
(iii) If $D \geq 4$, then $\theta_1 \geq  \frac{a_1 + \sqrt{a_1^2 + 4k}}{2}$.
\end{theorem}
\pf
Let $x \in V(\Gamma)$.  Then $\widehat{\Delta(x)}$ contains the equitable partition $\Pi=\{ \{ x\}, \Gamma(x)\}$ with quotient matrix
$Q(\Pi) =\left[
\begin{array}{cc}
 0  & k \\
 1 & a_1

\end{array} \right].$

As $Q(\Pi)$ has eigenvalues $\frac{a_1 \pm \sqrt{a_1^2 + 4k}}{2}$, it follows that $\frac{a_1 \pm \sqrt{a_1^2 + 4k}}{2}$ are eigenvalues of $\widehat{\Delta(x)}$.\\
\\
(i): By Theorem \ref{0}, one obtains $\theta_D \leq \frac{a_1 - \sqrt{a_1^2 + 4k}}{2}$. If $\theta = \frac{a_1 - \sqrt{a_1^2 + 4k}}{2}$ is an eigenvalue of $\Gamma$ then the standard sequence for $\theta$ satisfies $u_2(\theta) = 0$, and hence $\theta \neq \theta_D$ as the standard sequence for $\theta_D$ has exactly $D$ sign changes (\cite[Proposition 4.1.1]{bcn}).
This shows (i).
\\
(ii): This was shown by Koolen and Park \cite{Koolen Park}, but we give a proof for the convenience of the reader.
 The numbers $\frac{a_1 + \sqrt{a_1^2 + 4k}}{2}$ and $a_3$ are eigenvalues of the induced subgraph on $\{x\} \cup \Gamma(x) \cup \Gamma_3(x)$. Hence (ii) follows by Theorem \ref{0}.\\
(iii) Let $x$ and $y$ be two vertices at distance four. Then the induced subgraph of $\Gamma$ on $\{x, y\} \cup \Gamma(x) \cup \Gamma(y)$ is the disjoint union of $\widehat{\Delta(x)}$ and $\widehat{\Delta(y)}$. This graph has  $\frac{a_1 + \sqrt{a_1^2 + 4k}}{2}$ as an eigenvalue with multiplicity at least two, so by Theorem \ref{0}, it follows that $\theta_1 \geq \frac{a_1 + \sqrt{a_1^2 + 4k}}{2}$.
This finishes the proof of the theorem.
\epf

For diameter three we need some more specific information on the eigenvalues.
\begin{prop}\label{propb1}
Let $\Gamma$ be a distance-regular graph with diameter three and distinct eigenvalues $k = \theta_0 > \theta_1 > \theta_2 > \theta_3$. Then the following holds:\\
(i) $\theta_2 $ lies between $-1$ and $a_3 - b_2$;
\\
(ii) $\theta_1 > a_3 -b_2 > \theta_3$;
\\
(iii) $\theta_1 > a_1 -c_2 +1 > \theta_3$.
\end{prop}
\pf
Recall that
$$T =\left[
\begin{array}{ccc}
 -1 & b_1 & 0 \\
 1 & k- b_1 -c_2 & b_2 \\
0 & c_2 & a_3 - b_2

\end{array} \right]$$
and $\theta_1, \theta_2$ and $\theta_3$ are the three eigenvalues of $T$ . By Theorem \ref{0}, $\theta_2$ lies between $-1$ and $a_3 -b_2$. This shows (i).
\\
(ii): By Theorem \ref{0}, we have $\theta_1 \geq a_3 - b_2 \geq \theta_3$, but any eigenvector for $\theta_1$ of $T$ has only positive or only negative entries, and no zeroes (as $T$ is tridiagonal and the proof of \cite[Proposition 4.1.1]{bcn} still applies). This shows $\theta_1 \neq a_3-b_2$.
For any eigenvector for $\theta_3$ of $T$ has exactly two sign changes, so again 0 cannot occur in the eigenvector. This shows $\theta_3 \neq a_3 -b_2$.
\\
(iii):  Note that $\theta_1 + \theta_2 + \theta_3 = a_1 + a_2 + a_3 -k$. We need to consider two cases: $a_3 - b_2 \leq -1$ and $a_3 -b_2 > -1$.
First we consider the case  $a_3 - b_2 \leq -1$. Then by Theorem \ref{0}, $\theta_3 < a_3-b_2 \leq \theta_2 \leq -1 < \theta_1$, where the first inequality holds by (ii) and the last inequality holds by $\theta_1>0$.
We obtain
$$\theta_1 = a_1 +a_2 +a_3 - k - \theta_2 - \theta_3 > a_1 + a_2 + a_3 - k +1 + b_2 -a_3 = a_1 -c_2 + 1,$$ and $$\theta_3 = a_1 +a_2 +a_3 - k - \theta_1 - \theta_2
< a_1 + a_2 + a_3 - k +1 + b_2 -a_3 = a_1 -c_2 + 1.$$\\
Now we consider the case $a_3 - b_2 >-1$.
As any connected graph, with at least two vertices, has smallest eigenvalue at most $-1$ with equality if and only if the graph is complete (as an edge has eigenvalues $+1, -1$ and Theorem \ref{0}), $a_3-b_2>-1$ implies that $\theta_3 < -1 \leq \theta_2 \leq a_3 - b_2 < \theta_1$, where the last inequality holds by (ii). Then we obtain, in similar manner as above, $\theta_1 > a_1 - c_2 +1$ and $\theta_3 < a_1 -c_2 +1$. \epf

 For a connected graph $\Gamma$, the {\em  distance-$i$ graph} $\Omega$ has vertex set $V(\Gamma)$ and two vertices $x$ and $y$ are adjacent in $\Omega$ if they are at distance $i$ in $\Gamma$. Now we will extend  \cite[Proposition 4.2.17]{bcn} slightly.

\begin{prop}\label{-1}
Let $\Gamma$ be a distance-regular graph with diameter three and distinct eigenvalues $k = \theta_0 > \theta_1 > \theta_2 > \theta_3$.
Then the following are equivalent:
\\
\\
(i) $\theta_2 = -1$;\\
(ii) $\theta_2 = a_3 -b_2$;\\
(iii) $k+1 = c_3 + b_2$;\\
(iv) the distance-3 graph of $\Gamma$ is strongly regular.
\end{prop}
\pf By \cite[Proposition 4.2.17]{bcn}, (i), (iii) and (iv) are equivalent.
Also it is clear that (i) and (iii) together imply (ii). So we only need to show that (ii) implies (iii). Let $\theta_2 = a_3 - b_2$. As $\theta_1, \theta_2$ and $\theta_3$ are eigenvalues of $T$, hence $-b_2 c_2 (-1 - a_3 + b_2) = \ $det$(T - (a_3 - b_2)I) =0$. This shows that $-1-a_3+b_2 =0$ which is equivalent with (iii). This finishes the proof.
\epf

A distance-regular graph with diameter three is called {\em Shilla} if its second largest eigenvalue is equal to
 $\frac{a_1 + \sqrt{a_1^2 + 4k}}{2}$ (\cite{Koolen Park}).

Koolen and Park \cite{Koolen Park} gave the following characterization of Shilla distance-regular graphs.

 \begin{lemma}\label{shilla}
 Let $\Gamma$ be a distance-regular graph with diameter three,  valency $k$ and second largest eigenvalue $\theta_1$.
 Then the following are equivalent:
 \\
 (i) $\Gamma$ is Shilla, i.e., $\theta_1 = \frac{a_1 + \sqrt{a_1^2 + 4k}}{2}$;\\
 (ii) $k = a_3 (a_3 - a_1);$
 \\
 (iii) $\theta_1 = a_3$.
 \end{lemma}

\begin{prop}
Let $\Gamma$  be a distance-regular graph with $D$ at least three and distinct eigenvalues $k = \theta_0 > \theta_1 > \ldots > \theta_D$.
Then $\theta_1 = \frac{a_1 + \sqrt{a_1^2 + 4k}}{2}$ if and only if one of the following holds:
\\
(i) $D=3$ and $\Gamma$ is a Shilla distance-regular graph;
\\
(ii) $D=4$ and $\Gamma$ is an antipodal distance-regular graph.
\end{prop}
\pf
For $D=3$, this follows by its definition.\\
Let $D \geq 4$. Note that $\theta_1 = \frac{a_1 + \sqrt{a_1^2 + 4k}}{2}$ implies $u_2(\theta_1) = 0 > u_3(\theta_1) > \ldots > u_D(\theta_1)$ by \cite[p.130]{bcn}. By the Perron-Frobenius Theorem \cite[Theorem 3.1.1]{bcn} we obtain that $\theta_1$ is the largest eigenvalue of
$$\left[\begin{array}{cccc}
 a_3 & b_3 &  & {\bf 0}\\
 c_4 & a_4 & b_4  \\
 & \ddots & \ddots & \ddots \\
{\bf 0} &  & c_D & a_D

\end{array} \right]$$
as the $b_i$'s and $c_i$'s are positive integers.

Let $\Sigma(x)$ be the induced subgraph of $\Gamma$ on $\Gamma_3(x) \cup \ldots \cup \Gamma_D(x)$. Let $y \in \Gamma_4(x)$. Then $\widehat{\Delta(y)}$ is a subgraph of $\Sigma(x)$ and has largest eigenvalue $\frac{a_1 + \sqrt{a_1^2 +4k}}{2}$. Let $\Lambda$ be the connected component of $\Sigma(x)$ that contains $y$.  Then the following hold:\\
(a) The largest eigenvalue of $\Lambda$ is at most
$\frac{a_1 + \sqrt{a_1^2 +4k}}{2}$; \\
(b) $\Lambda$ contains $\widehat{\Delta(y)}$ as an induced subgraph; and\\
(c) $\widehat{\Delta(y)}$ has eigenvalue
$\frac{a_1 + \sqrt{a_1^2 +4k}}{2}$. \\
This means that by the Perron-Frobenius Theorem  $\widehat{\Delta(y)}$ must be equal to $\Lambda$ and hence is a connected component of $\Sigma(x)$.
This implies that, $D=4, c_4 = k, b_3 = 1$ and $a_1 = a_3$, as $y$ is the only vertex with degree $k$ in $\widehat{\Delta(y)}$. This concludes the proof.
\epf

\noindent
{\bf Remark.} An antipodal distance-regular graph with diameter four has second largest eigenvalue $\frac{a_1 + \sqrt{a_1^2 +4k}}{2}$ and second smallest eigenvalue $\frac{a_1 - \sqrt{a_1^2 +4k}}{2}$, as these are the eigenvalues of its folded graph (which is strongly regular).
\\
\\
The following result is the main result in this section.
\begin{theorem}\label{upp}
Let $\Gamma$ be a distance-regular graph with diameter $D$ at least two, and distinct eigenvalues $k = \theta_0 > \theta_1 > \ldots > \theta_D$.
Then $$(\theta_1 +1)(\theta_D +1) \leq -b_1$$ holds with equality if and only if the diameter $D$ is equal to two.
\end{theorem}

\pf
Let $\Gamma$ be a distance-regular graph with diameter $D$ at least two. \\
For $D=2$, it is well-known that $\theta_1=\frac{a_1-c_2+\sqrt{(a_1-c_2)^2+4(k-c_2)}}{2}$ and $\theta_D=\frac{a_1-c_2-\sqrt{(a_1-c_2)^2+4(k-c_2)}}{2}$ (see for example \cite[p.220]{godsil-combin}).
So this shows that the theorem holds for $D$ is equal to 2.

For $D=3$, $\theta_1+1, \theta_2+1$ and $\theta_3+1$ are the eigenvalues of the matrix
\begin{equation}\label{inter}
T+I = \left[
\begin{array}{ccc}
 0 & b_1 & 0 \\
 1 & k+1 - b_1 -c_2 & b_2 \\
0 & c_2 & k+1  - b_2-c_3
\end{array} \right].
\end{equation}
If $\theta_2=-1$, then $(\theta_1+1)(\theta_D+1)= -b_1 -b_2c_2 < -b_1$, by looking at the coefficient of the linear term in the characteristic polynomial of $T + I$.
If $\theta_2\neq-1$, then $(\theta_1+1)(\theta_2+1)(\theta_3+1)=-b_1(k-b_2-c_3+1)$ and $\mid \theta_2+1 \mid$  $\leq$  $\mid k-b_2-c_3+1 \mid$ by Theorem \ref{0}. But the equality case can not hold by Proposition \ref{-1}. Therefore the inequality follows in this case.

For $D\geq 4$, the inequality follows immediately from Theorem \ref{ineq1}.
\epf

\noindent
{\bf Remark.} An antipodal distance-regular $r$-cover with diameter three satisfies $(\theta_1 +1)(\theta_3+1) =
-b_1\frac{r}{r-1}$.  De Caen, Mathon, and Moorhouse \cite{CMM95} constructed
distance-regular antipodal $2^{2t-1}$-cover of the complete
graph $K_{2^{2t}}$, i.e., with intersection array
$\{2^{2t}-1,2^{2t}-2,1;1,2,2^{2t}-1\}$.  This shows that we can not improve the bound for $D=3$ in the above theorem. It is likely that the above inequality can be improved for larger diameter.\\

\begin{cor} Let $\Gamma$ be a distance-regular graph with diameter $D$ at least three and distinct eigenvalues $k=\theta_0 > \theta_1 > \ldots > \theta_D$.
Let $x$ be a vertex of $\Gamma$ and let $\eta_1 = a_1 \geq \eta_2 \geq \ldots \geq \eta_k$ be the eigenvalues of the local graph $\Delta(x)$. Then
$\theta_1 > \eta_2 > \eta_k > \theta_D$.
\end{cor}
\pf Terwilliger \cite[cf. Theorem 4.4.3]{bcn} showed that $\eta_2 \leq -1 - \frac{b_1}{\theta_D +1}$ and   $\eta_k \geq -1 - \frac{b_1}{\theta_1 +1}$. Now the corollary immediately follows from Theorem \ref{upp}, as the local graph is not complete ($D \geq 3$).
\epf

\section{Distance-regular graphs with fixed second largest eigenvalue}
In this section we look at  distance-regular graphs with diameter at least three and fixed second largest eigenvalue $\theta_1$. \\
First we will look at the case when  $\theta_1 \leq 1$ and show that then the graph must be a $K_{n,n}$ in which a perfect matching is removed.

Moreover, we will show that, for a fixed integer  $m$ at least two, there are only finitely many distance-regular graphs with diameter at least three and valency at least three such that its second largest eigenvalue is between 1 and $m$. We conclude this section by determining all the distance-regular graphs with diameter at least three, valency at least three and second largest eigenvalue at most two.  Note that the situation for distance-regular graphs with diameter at least three is completely different from the situation for strongly regular graphs, as for every positive integer $m$ there are infinitely many connected strongly regular graphs with second largest eigenvalue $m$.

First we look at the case   $\theta_1\leq 1$:

\begin{theorem}
Let $\Gamma$ be a distance-regular graph with valency $k$ at least three, diameter $D$ at least three and distinct eigenvalues $k= \theta_0 > \theta_1> \ldots> \theta_D$  satisfying $\theta_1 \leq 1$. Then $\Gamma$ has diameter three and is the graph $K_{k+1,k+1}$ in which a perfect matching is removed.
\end{theorem}
\pf  As $k \geq 3$,  it follows that $\frac{a_1+\sqrt{a_1 ^2+4k}}{2}\geq \sqrt{3} > 1$. So $\theta_1 <   \frac{a_1+\sqrt{a_1 ^2+4k}}{2}$, and hence $D = 3$ and
$1\geq\theta_1\geq a_3$ hold by Theorem \ref{ineq1}. If $a_3 =1$ holds, then $\theta_1 = 1$, which is impossible by  Lemma \ref{shilla}. Hence we obtain $a_3=0$. Then, $\theta_0\theta_1\theta_2\theta_3=b_2 k^2$, as $\theta_0, \theta_1, \theta_2$ and $\theta_3$ are the eigenvalues of
$$\left[
\begin{array}{cccc}
 0 & k & 0 & 0 \\
 c_1 & a_1 & b_1&  0   \\
 0& c_2 & a_2 & b_2 \\
 0 & 0 & c_3 & a_3
\end{array} \right] \ \  (\cite[\mbox{p.} 129]{bcn}).$$ Now
$\mid\theta_2\mid\leq b_2$ (Proposition \ref{propb1})  implies that $\theta_3\leq -k$ (as
$\theta_1\leq1)$. This show that $\Gamma$ is bipartite and
$\theta_2=-b_2$. Now, it follows from Proposition \ref{-1} that $b_2=1$ holds. This shows the theorem. \epf

Now we consider the case when $\theta_1$ is more then 1, but at most a fixed positive integer at least two.

\begin{theorem}\label{class3}
Let $m \geq 2$ be an integer. Then there are finitely many distance-regular graphs with valency $k$ at least three, diameter $D$ at least three and distinct eigenvalues $k= \theta_0 > \theta_1>  \ldots> \theta_D$ such that $1<\theta_1\leq m$.
\end{theorem}
\pf
If $\theta_1\geq \frac{a_1+\sqrt{a_1^2+4k}}{2}$, then
$m \geq \frac{a_1+\sqrt{a_1^2+4k}}{2} \geq \sqrt{k}$, so $k \leq m^2$. Bang et al. \cite{BDKM} showed that there are finitely many distance-regular graphs with valency $k$ such that $3 \leq k \leq m^2$.
So we may assume that $\theta_1 < \frac{a_1+\sqrt{a_1^2+4k}}{2}$ and hence $D=3$ and
$m\geq\theta_1\geq a_3$ hold by Theorem \ref{ineq1}.
First we consider the case that $\Gamma$ is not bipartite. Then, by Theorem \ref{upp}, we obtain $(\theta_3+1)\leq \frac{-b_1}{m+1}$. Next, we show  that, $b_1 \geq \frac{k+1}{3}$ holds. Let $y$ and $z$ be two vertices at distance two, and let $x$ be a common neighbor of $y$ and $z$. Then $
c_2 - 1 \geq |\Gamma(x) \cap \Gamma(y) \cap \Gamma(z)| \geq
2(a_1 + 1) - k$. As $ D = 3$, we have $c_2 \leq b_1$ (Lemma \ref{pre} (iii)). It follows that $a_1  \leq \frac{2(k-2)}{3}$ and hence $b_1 \geq \frac{k+1}{3}$.
This implies $\theta_3 < \frac{-k-1}{3(m+1)}$.
Recal that $m_i$ is the multiplicity of $\theta_i$ $(i=0,1,2,3)$.
Then $0 \leq \sum_{x \in V(\Gamma)}(A^3)_{xx} = \sum_{i =0}^3 m_i \theta_i^3$ and,
as $\theta_3 < \frac{-k-1}{3(m+1)}$, $\theta_2 < \theta_1 \leq m$ and $(m_3-1)(m_3 +2)  \geq 2k$ (as $\Gamma$ is not bipartite (cf. \cite[Theorem 5.3.2]{bcn})) this implies
$0 < k^3 + (v-1)m^3 + (\sqrt{2k}-1)(\frac{-k-1}{3(m+1)})^3$ (where $v = |V(\Gamma)|$). This means that $k$ is bounded above by a polynomial in $m$. So this shows the theorem in case  when $\Gamma$ is not bipartite.
The remaining case is when $\Gamma$ is bipartite and $D=3$. In this case $\theta_1 = \sqrt{b_2}$ and hence  $2 \leq b_2 \leq m^2$. Now for a fixed number $b_2\geq 2$, we know that $k-b_2$ divides $k(k-1)$ (as $k_2$ is an integer) and this implies $k-b_2$ divides $b_2(b_2-1)$. Since $b_2\geq 2$, we obtain $k\leq(b_2)^2$. This completes the proof of the theorem.
\epf

Now we will determine the distance-regular graphs with dameter at least three, valency at least three and whose second largest eigenvalue is more then one but at most two.

\begin{theorem}
Any distance-regular graph with valency $k$ at least three, diameter $D$ at least three
and distinct eigenvalues $k =\theta_0 > \theta_1>  \ldots> \theta_D$  such that $1<\theta_1\leq2$ has one of the 23
intersection arrays listed in Table 1 below.
\end{theorem}
\pf  The (intersection arrays of) distance-regular graphs with
valency three and four were classified by Biggs et al. \cite{Biggs}
and Brouwer and Koolen \cite{KB}, respectively. By checking them, we
obtain the intersection arrays 1--3,13--19 and 23. So we may assume
that $k \geq 5$ in the rest of the proof. This means that $\theta_1<
\frac{a_1+\sqrt{a_1^2+4k}}{2}$ as $\theta_1 \leq 2$ and hence $D=3$
and $2 \geq \theta_1 > a_3$ by Theorem \ref{ineq1} (iii) and Lemma
\ref{shilla}.
So this implies that $a_3\in\{0,1\}$.  \\
Now we consider two cases, namely $b_2 = a_3$ and $b_2 \neq a_3$.
Suppose first that $b_2=a_3$. In this case we will show that $k
\leq 12$. Note that $b_2 = a_3$ implies  $b_2=a_3=1$ (as $b_2 \geq
1$) and as $\theta_1, \theta_2$ and $\theta_3$ are eigenvalues of
$T$ we find $\theta_1\theta_2\theta_3=c_2>0$  and  $-1\leq\theta_2
< 0$. Because $k_3 = \frac{kb_1}{c_2(k-1)}$ is an integer, it follows $k-1$ divides $kb_1$, and hence $b_1 = k-1$. Now by
Theorem \ref{upp}, we have $\theta_3\leq\frac{-k-2}{3}$. As
$-c_2 = a_1 + a_2 + a_3 -k=\theta_1+\theta_2+\theta_3 \leq 2+\frac{-k-2}{3}$, it
follows that $c_2\geq \frac{k-4}{3}$. This gives us
$k_3=\frac{k}{c_2}\leq \frac{3k}{k-4}$. If $k \geq 11$, then $k_3
\leq 4$ and, as $\Gamma$ is not antipodal, we see $k \leq k_3(k_3-1)
\leq 12$, by \cite[Proposition 5.6.1]{bcn}. So, if $b_2 = a_3$, then
$k \leq 12$ holds. Now we consider the case  that $b_2 \neq a_3$. This
implies that $b_2 > a_3$, as $a_3 \leq 1$. We will show that $k \leq
25$ holds in this case. Since $-1 \geq \theta_2 \geq a_3 - b_2$ (Proposition 3.2(i)),
$1 < \theta \leq2$, and
$\theta_1\theta_2\theta_3=k(b_2-a_3)+a_3c_2$ hold, we find
$\theta_3 \leq-k/2$.
We will now show that, $m_3 \geq k/2$ or $\Gamma$ is bipartite, hold. To do so, we assume $m_3 < k/2$. Then  \cite[Theorem
4.4.4]{bcn} implies that $\frac{b_1}{\theta_3+1}$ is an integer.
Since $-1\geq\frac{b_1}{1+\theta_3}\geq \frac{2-2k}{k-2}$ and $k\geq
5$ holds, we obtain $\frac{b_1}{1+\theta_3}\in \{-1, -2\}$. In the case when
$\frac{b_1}{1+\theta_3}=-1$, then $\Gamma$ is bipartite by
\cite[Proposition
 4.4.7]{bcn}, as
$u_2(\theta_3)=1$.
 In the case when $\frac{b_1}{1+\theta_3}=-2$, we obtain
$a_1\leq 1$, as $\theta_3 \leq -k/2$. By \cite[Theorem 4.4.4]{bcn}, we see that $-1 - \frac{b_1}{1+\theta_3} = 1$ is an eigenvalue of any local graph $\Delta(x)$ with multiplicity more then $k/2 $. But this implies that $a_1 =1$, but this gives a contradiction as then $1$ has exactly multiplicity $k/2$ in $\Delta(x)$.  So we have shown that $\Gamma$ is bipartite or $m_3 \geq k/2$.

Let us first consider the case that $\Gamma$
is not bipartite. Then the number of vertices is bounded above by $2k^2$
as $k_2 \leq k(k-1)$ and $k_3 = \frac{b_2}{c_3} k_2 < k_2$.  By
considering the trace of $A^3$, we find $k^3 + 2k^2 2^3 +
\frac{k}{2} (\frac{-k}{2})^3 > 0$. This implies $k \leq 25$. Now if
$\Gamma$ is bipartite, then $\theta_1 = \sqrt{b_2}$ and as
$1<\theta_1  \leq 2$, we have $b_2 \in \{2,3,4\}$. Now  $k_2 =
\frac{k(k-1)}{c_2} = \frac{k(k-1)}{k-b_2}$ must be an integer and
this implies $k \leq 16.$ In conclusion, we obtain that if $k \geq 5$,
then $D =3$ and $k \leq 25$. By computer checking, we find that the arrays
4-12 and 20-22 in Table 1 are the only possible intersection arrays.
This completes the proof. \epf

\hspace*{-0.7cm}
\begin{tabular}{cccllc}\hline
 No. & $v$ & $D$ & intersection array & spectrum & $\begin{matrix} \mbox{number \ of} \\ \mbox{non-isomorphic  \ graphs} \end{matrix}$ \\ \hline
1&14&3&\{3,2,2;1,1,3\}& $\pm(3^1 \ \sqrt2 ^6)$ & 1 \\
2&14&3&\{4,3,2;1,2,4\}&$\pm(4^1 \ \sqrt2 ^6)$ & 1 \\
3&26&3&\{4,3,3;1,1,4\}& $\pm(4^1 \ \sqrt3 ^{12})$ & 1 \\
4&22&3&\{5,4,3;1,2,5\}& $\pm(5^1 \ \sqrt3 ^{10})$ & 1 \\
5&42&3&\{5,4,4;1,1,5\}& $\pm(5^1 \ 2^{20})$ & 1 \\
6&22&3&\{6,5,3;1,3,6\}& $\pm(6^1 \ \sqrt3 ^{10})$ & 1 \\
7&32&3&\{6,5,4;1,2,6\}& $\pm(6^1 \ 2^{15})$ & 3 \\
8&30&3&\{7,6,4;1,3,7\}& $\pm(7^1 \ 2^{14})$ & 5 \\
9&30&3&\{8,7,4;1,4,8\}& $\pm(8^1 \ 2^{14})$ & 5 \\
10&26&3&\{9,8,3;1,6,9\}& $\pm(9^1 \ \sqrt3 ^{12})$ & 1 \\
11&32&3&\{10,9,4;1,6,10\}& $\pm(10^1 \ 2^{15})$ & 3 \\
12&42&3&\{16,15,4;1,12,16\}& $\pm(16^1 \ 2^{20})$ & 1 \\
13&18&4&\{3,2,2,1;1,1,2,3\}& $\pm(3^1 \ \sqrt3 ^6) \ 0^4$ & 1 \\
14&30&4&\{3,2,2,2;1,1,1,3\}& $\pm(3^1 \ 2^9) \ 0^{10}$ & 1 \\
15&16&4&\{4,3,2,1;1,2,3,4\}& $\pm(4^1 \ 2^4) \ 0^6$ & 1 \\
16&32&4&\{4,3,3,1;1,1,3,4\}& $\pm(4^1 \ 2^{12}) \ 0^6$ & 1 \\
17&20&5&\{3,2,2,1,1;1,1,2,2,3\}& $\pm(3^1 \ 2^4 \ 1^5)$ & 1 \\
18&15&3&\{4,2,1;1,1,4\}& $4^1 \ 2^5 \ (-1)^4 \ (-2)^5$ & 1 \\
19&35&3&\{4,3,3;1,1,2\}& $4^1 \ 2^{14} \ (-1)^{14} \ (-3)^6$ & 1 \\
20&36&3&\{5,4,2;1,1,4\}& $5^1 \ 2^{16} \ (-1)^{10} \ (-3)^9$ & 1 \\
21&42&3&\{6,5,1;1,1,6\}& $6^1 \ 2^{21} \ (-1)^{6} \ (-3)^{14}$ & 1 \\
22&27&3&\{8,6,1;1,3,8\}& $8^1 \ 2^{12} \ (-1)^{8} \ (-4)^6$ & 2 \\
23&28&4&\{3,2,2,1;1,1,1,2\}& $3^1 \ 2^8 \ (\pm\sqrt2-1)^6 \ (-1)^7 $ & 1 \\
\hline
& & & & & \\
\multicolumn{6}{c}{Table 1: The intersection arrays of the distance-regular graphs with}\\
\multicolumn{6}{c}{  diameter at least three and $1 < \theta_1 \leq 2$  \  \ \ \ }
 \end{tabular} \\
\newpage
\noindent
{\textit{Descriptions of these graphs}}
Note that a bipartite distance-regular graph $\Gamma$ with intersection array $\{k, k-1, k-c_2; 1, c_2, k\}$ is the point-block incidence graph of a square $2-(\frac{n}{2}, k, c_2)-$design where $n$ is the number of vertices of $\Gamma$, and vice versa  the point-block incidence graph of a square $2-(v, k, \lambda)-$design is distance-regular with intersection array $\{k, k-1, k-\lambda; 1, \lambda, k\}$. The first 12 graphs in the table are bipartite with diameter three and fall into this class.
1. the Heawood graph; 13. the Pappus graph; 14. Tutte's 8-cage; 15. the 4-cube; 16. the incidence
graph of the $AG(2,4)$ minus a parallel class; 17. the Desargues graph;
18. the line graph of Petersen graph; 19. the Odd graph $O_4$; 20. the Sylvester
graph; 21. the subgraph graph induced on the second subconstituent of the Hoffman-Singleton graph; 22. there are two
non-isomorphic graphs; one of them is the $GQ(2,4)$ minus a spread; 23. the
Coxeter graph.
\\
\\
{\bf Acknowledgements}
We would like he anonymous referee for their comments as they greatly improved the paper. Also we would like to thank Sejeong Bang for her careful reading of the paper.

\end{document}